\documentclass[notitlepage, 10pt]{article}
\usepackage{hyperref, amsmath, amsfonts, enumitem, amsthm, xcolor, parskip, stmaryrd}
\usepackage{amssymb} 
\usepackage[nottoc]{tocbibind}
\usepackage{tikz-cd}
\usepackage{thmtools, cleveref}
\usepackage{subfiles}

\hypersetup{             
    colorlinks=true,                
    breaklinks=true,                
    urlcolor= black,                
    linkcolor= blue,                     
    bookmarksopen=false,
    filecolor=black,
    citecolor=blue,
    linkbordercolor=blue
}
\setlist[itemize]{itemsep=0pt, parsep=0pt}
\setlist[enumerate]{itemsep=0pt, parsep=0pt}

\newtheoremstyle{theo}
  {0.5em} 
  {\topsep} 
  {\itshape} 
  {} 
  {\bfseries} 
  {.} 
  {0.2em} 
  {} 

\newtheoremstyle{defin}
  {0.5em} 
  {\topsep} 
  {} 
  {} 
  {\bfseries} 
  {.} 
  {0.2em} 
  {} 

\newtheoremstyle{rema}
  {0.5em} 
  {\topsep} 
  {} 
  {} 
  {\scshape} 
  {.} 
  {0.2em} 
  {} 

\theoremstyle{defin}
\newtheorem{defn}{Definition}
\newtheorem{eg}{Example}[section]
\theoremstyle{theo}
\newtheorem{thm}[defn]{Theorem}

\newtheorem{lem}{Lemma}[defn]
\newtheorem{prop}[defn]{Proposition}
\newtheorem{cor}{Corollary}[defn]
\theoremstyle{rema}

\newtheorem{rem}{Remark}[section]
\newcommand{\R}{\mathbb{R}}

\newcommand{\N}{\mathbb{N}}
\newcommand{\Z}{\mathbb{Z}}

\newcommand{\bB}{\mathbb{B}}
\newcommand{\cC}{\mathcal{C}}
\newcommand{\cM}{\mathcal{M}}
\newcommand{\cD}{\mathcal{D}}
\newcommand{\cO}{\mathcal{O}}

\newcommand{\cv}{\mathrm{conv}}

\title{Finitely Generated Congruences in Semirings and Canonical Positive Models}
\author{Snehinh Sen\footnote{\textsc{Mathematical Sciences Institute\\ The Australian National University\\ Canberra, ACT, Australia - 2601}\\
E-mail: snehinh.sen@anu.edu.au}}
\setlist[enumerate, 1]{label={(\roman*)}}

\begin{document}

\maketitle

\begin{abstract}
In this paper, we explore semirings in which all congruences are finitely generated. Such semirings are dubbed \textit{Congruence Noetherian}. After developing sufficient background and examples, we focus on the canonical positive models of a real order and show that this obvious choice, though not finitely generated as an $\N$-module, is both Congruence Noetherian and flat over $\N$.
\end{abstract}

\textit{Keywords:} Semirings, Noetherian, Hilbert Basis Theorem, Flatness, Congruences, Positive Models

\textit{MSC 2020:} Primary: 16Y60 , 13F20.

\section{Introduction}\label{s1}

A \textit{semiring} is, essentially, a ring in which some elements might not have additive inverses. We still assume the existence of the respective operational identities $0, 1$ as well as that $0\cdot a = 0$ for each element $a$ of the semiring. Finding interesting general properties for semirings is difficult because of their versatility. Nevertheless, one is often tempted to consider the most obvious generalizations of the usual definitions from the world of rings. For example, one may define a semiring to be \textit{Noetherian} if every ideal is finitely generated. This notion has been explored, though not extensively (see, for instance, \cite{abuhlail2019knoetheriankartiniansemirings}, \cite{allencohen}, \cite{allenideal2006}, \cite{allenhilb1980}). Interestingly, for semirings, there is another obvious generalization of Noetherianness. More precisely, we can look at semirings for which every congruence is finitely generated. Such semirings will be called \textit{congruence Noetherian} or \textit{c-Noetherian}. A version of this has appeared in \cite{JooMincheva2017} for additively idempotent semirings. Similarly, a semiring is said to be \textit{congruence principal} or \textit{c-principal} if every congruence is generated by a single relation.

The first part of this article studies the properties of c-Noetherian semirings and provides ample examples. One of the main results that we prove is the following.

\begin{restatable}{thm}{thilb}[Hilbert Basis Criterion]\label{thm1}
        Let $R$ be a semiring. The following are equivalent.
    \begin{enumerate}
        \item $R[X]$ is c-Noetherian.
        \item Every finitely generated algebra over $R$ is c-Noetherian.
        \item Every finitely generated algebra over $R$ is finitely presented.
        \item $R$ is a Noetherian ring. 
    \end{enumerate} 
\end{restatable}

Even though attempts to generalize Hilbert Basis Theorem exist in the literature on semirings (see, for instance, \cite{chen2026arithmeticsemiringsiideals}, \cite{allenhilb1980}), to the best of our knowledge, the above result and the constructions presented here have not been explored. 

In the second part of this article, we shift our focus to subgroups and subrings of real numbers and their positive models. In this setting, a positive model is a corresponding object defined over $\N$ (instead of $\Z$) that recovers the initial object after base change. The most natural positive model is the set of all non-negative elements $S=R_{\geq 0}$ of a subgroup $R$. As this model is independent of choice of basis, it will be called the \textit{canonical positive model} of $R$. After verifying functoriality, we prove a few fascinating properties of $S$. Recall that a \textit{k-ideal} of a semiring $R$ is a cancellative ideal. We say that a semiring is \textit{k-Noetherian} if its k-ideals satisfy the ascending chain condition. 

\begin{restatable}{thm}{ttarget}[Positive Models of Subrings of $\R$]\label{t:target}
    Let $R$ be a subgroup of $\R$ and $S=R_{\geq 0}$. Then $S$ is flat over $\N$. Furthermore, if $R$ is a subring, then
    \begin{enumerate}
        \item $S$ is c-Noetherian if and only if $S$ is k-Noetherian, which is if and only if $R$ is Noetherian.
        \item $S$ is c-principal if and only if $R$ is a principal ideal domain.
        \item If $R$ is a finitely generated $\Z$-module (that is, a \textit{real order}), then every non-trivial quotient is finite. 
    \end{enumerate}
\end{restatable}

These are, in some sense, the best that one could hope for (as, for instance, $S$ is not finitely generated as an $\N$-module, nor is $S$ free, even if $R$ is free).

Regarding the structure of the article, Section \ref{s2} describes some preliminaries and known results from the theory of semirings. Section \ref{s3} studies c-principal and c-Noetherian semirings. Section \ref{s4} gives a somewhat universal example of a semiring which is not c-Noetherian and we prove Theorem \ref{thm1}. We shift our focus to canonical positive models to prove their finiteness properties in Section \ref{s7} and their flatness in Section \ref{s8}.

\section{Basics} \label{s2}

Let $R$ be a semiring and $\cC\subseteq R\times R$ be an equivalence relation denoted by $\sim$, that is, we write $a\sim b$ if $(a,b)\in \cC$. We say that $\cC$ is a \textit{congruence} if $a\sim b$, $c\sim d$ implies $a+c \sim c+d$ and $ac \sim bd$. As $\cC$ is also reflexive, this includes the cases of $ar\sim br$ and $a+r\sim b+r$ for any $r\in R$. In other words, it is a subsemiring of $R\times R$ which is also an equivalence relation. 

It is clear that the intersection of congruences is a congruence, so we may talk about the congruence generated by a subset $B\subseteq R\times R$. This is denoted by $\langle a\sim b : (a,b)\in B \rangle$ If $B$ is empty, this is simply the diagonal congruence. We say that a congruence is \textit{finitely generated} if it can be generated by a finite set.

As with rings, all examples of congruences come from quotients.

\begin{eg}
    Let $f : R \to S$ be a surjective semiring homomorphism. Then the \textit{congruence kernel of $f$}, denoted by $\cC_f$ is defined as $\{(a,b) : f(a)=f(b)\}$. It is easy to verify that this is indeed a congruence. Conversely, given a congruence $\cC$ of $R$, the equivalence classes $R/\cC$ of $\cC$ form a semiring under the operations induced from $R$. Furthermore, $\cC_{\pi : R \to R/\cC} = \cC$.
\end{eg}

The following proposition allows us to easily work with congruences with a prescribed generating set. This is the semiring version of the Mal'tsev Congruence Generation Theorem (see, for example, Theorem 6.17 of \cite{FreeseMcKenzieMcNultyTaylor2022}).

\begin{prop}
    \label{p:definingconggen}
    Let $A$ be a semiring and $\cC$ be a congruence on $A$ generated by pairs $(a_i, b_i)\in A^2$. Define an elementary relation as a relation of the form $ua_i+vb_i+c\sim ub_i+va_i+c$ with $u,v,c\in A$ and $i\in I$. Then given any $(a,b)\in \cC$, there exists an integer $n\geq 1$ and $r_0,\ldots, r_n \in A$ such that $r_0=a, r_1=b$ and $(r_{k-1}, r_{k})\in \cC$ is an elementary relation for each $k=1,\ldots, n$.
\end{prop}

\begin{proof}
    Define $\cC'$ to be the set of pairs $(a,b)\in A^2$ for which such a collection of elementary relations exist. It is clear that $\cC'\subseteq \cC$ and that $(a_i,b_i)\in \cC'$ for each $i\in I$. It suffices to prove that $\cC'$ is a congruence. That it is an equivalence relation is immediate. Furthermore, if $a = r_0\sim r_1\sim \ldots \sim r_n = b$ and $a' = r'_0\sim r_1'\sim \ldots \sim r'_{n'} = b'$ are two chains of elementary relations, then we observe that
    $$aa' = ar'_0 \sim ar'_1 \sim \ldots ar'_{n'} = r_0b' \sim r_1b' \sim \ldots \sim r_nb' =  bb'$$
    and 
    $$a+a' = a+r'_0 \sim a+r_1'\sim \ldots \sim a + r'_{n'} = r_0 +b' \sim r_1 + b'\sim \ldots r_n+b'=b+b'$$
    are chains of elementary relations, proving that $\cC'$ is a congruence.
\end{proof}

We now recall the definition of finitely generated and finitely presented algebras over a semiring and prove a general result. 

\begin{defn}\label{df2}
    Let $R$ be a semiring and $S$ be an $R$-algebra.
    \begin{enumerate}
        \item We say $S$ is \textit{finitely generated as an $R$-algebra} if there is a surjective map from a polynomial ring $R[X_1,\ldots,X_n]$ of finitely many variables to $S$.
        \item We say $S$ is \textit{finitely presented as an $R$-algebra} if there is a surjective map as described above with a congruence kernel finitely generated as a congruence.
    \end{enumerate}
\end{defn}

Just as in the case of rings, we have the following invariance result.

\begin{prop}
    \label{p3} Suppose $R$ is a semiring, $T$ is a finitely presented $R$-algebra,and $f : R[Y_1,Y_2,\ldots,Y_k] \to T$ is a surjective map. Then $C_f$ is finitely generated. 
\end{prop}

\begin{proof}
    Let $R[X_1,\ldots X_l]/C'\xrightarrow{\sim} T$ be a finite presentation of $T$ induced by a map $g$. Let $x_i\in T, y_j\in T$ be the images of $X_i$ and $Y_j$ with respect to $g$ and $f$ respectively. Let $h_j$ be some lifts of $y_j$'s in $R[\underline{X}] := R[X_1,\ldots,X_l]$ and $g_i$ be some lifts of $x_i$ in $R[\underline{Y}] :=R[Y_1,\ldots,Y_k]$. Consider the commutative diagram:

        \begin{equation}
            \begin{tikzcd}
                       & R[\underline{X}] \arrow[rd, twoheadrightarrow, "\cC'"]& \\
            R[\underline{X}, \underline{Y}] \arrow[ru, twoheadrightarrow, "\langle Y_j\sim h_j\rangle "] \arrow[rd, twoheadrightarrow, "\langle X_i\sim g_i \rangle"] & & T \\
                       & R[\underline{Y}] \arrow[ru, twoheadrightarrow, "\cC_f"]&       
            \end{tikzcd}
        \end{equation}
        
    Each arrow is marked with the corresponding kernel congruence. Then the kernel of the composition, say $\cC''$, is generated by $\cC'$ and $Y_j\sim h_j$, and is thus finitely generated. Thus, $\cC_f$, being the image of $\cC'$ under the map $R[\underline{X}, \underline{Y}]\twoheadrightarrow R[\underline{X}]$, is finitely generated.
\end{proof}

We briefly recall a classification result about semirings. Let $\bB$ be the Boolean semifield, defined as the quotient $\N/\langle 1 \sim 2\rangle$.

\begin{thm}[Borger -- Grinberg Theorem]\label{thBJ}
    Let $R$ be a semiring. Then the following are equivalent.
    \begin{enumerate}
        \item $R$ is a ring.
        \item $R\otimes_\N \bB = 0$
        \item $R$ does not admit $\bB$ as a quotient.
    \end{enumerate}
\end{thm}

\begin{proof}
    Proposition 7.9 of \cite{Borger_2015}.    
\end{proof}
Before concluding the section, we shift our focus to the flatness of modules over semirings and recall a few properties from \cite{borger_jun_semiring}. Let $R$ be a semiring and let $M$ be an $R$-module. Recall that $\mathrm{Hom}_R(M,-)$ has a left adjoint $-\otimes_R M$ in the category of $R$-modules that has many properties analogous to the usual tensor product over rings. For more details, see \cite{Borger2016}, \cite{borger_jun_semiring}, \cite{katsov2004}, and \cite{ToenVaquie2009}.

\begin{defn}
    \label{d:flat}
    We say $M$ is \textit{mono-flat} if $F_M = -\otimes_R M : \mathbf{Mod}_R \to \mathbf{Mod}_R$ preserves monomorphisms. We say $M$ is \textit{flat} if $F_M$ preserves finite limits.
\end{defn}

While all flat modules are mono-flat, the converse is not true generally (see, for instance, Example 3.7 in \cite{katsov2004}). The conditions are equivalent whenever $R$ is a ring or a Boolean algebra (see Theorem 3.2, \cite{katsov2004}). Theorem 3.5 of \cite{borger_jun_semiring} beautifully characterizes flat modules by giving multiple equivalent conditions. For the purpose of this talk, we will only need the following result.

\begin{cor}
    An $R$-module $M$ is flat if it is a directed union of its finite, free submodules.
\end{cor}

\begin{proof}
    Tensor product commutes with filtered colimits. A directed union is a filtered colimit.
\end{proof}

\begin{rem}
    Note that the converse is not true for a general $R$. In fact, it is not even true for an arbitrary ring. For instance, let $R$ be Dedekind Domain of class number greater than $1$ and $M$ be a non-principal ideal. Then $M$ is projective, hence flat over $R$, but it is definitely not the union of its finite free submodules. A simpler example comes from disconnected rings like $R=R_1\times R_2$ and $M = R_1\times 0$. In fact, For a modified converse, see Theorem 3.5 of \cite{borger_jun_semiring}. For convenience, we shall call such modules \textit{strongly flat}.
\end{rem}

Flatness is important for many reasons. It preserves a lot of necessary algebraic and geometric properties. At the same time, it is not that restrictive and hard to keep a track of like freeness. Therefore, it sits at the heart of the arithmetic and geometry of semirings. We refer the interested readers to \cite{Borger2016}, \cite{Borger_2015}, \cite{borger_jun_semiring}, \cite{ConnesConsani2011}, \cite{Zelich2026}, \cite{GualdiKuhrsMayoGarciaXarles2026}, \cite{katsov2004}, \cite{Lorscheid2017},  \cite{MaclaganSturmfels2015}, and \cite{ToenVaquie2009} for an overview of the theory with applications.

\section{Principal and Congruence Noetherian Semirings}\label{s3}

\begin{defn}\label{df1}
\begin{enumerate}
    \item We say that a congruence is \textit{principal} if it can be generated by a single relation of the form $a\sim b$.
    \item We say that a semiring is \textit{congruence principal} or \textit{c-principal} if every congruence is principal.
    \item We say that a semiring is \textit{congruence Noetherian} or \textit{c-Noetherian} if every congruence is finitely generated.
\end{enumerate}
\end{defn}

The following are examples of c-Noetherian and c-principal semirings.

\begin{eg}\label{e1}
    \begin{enumerate}
    \item A ring is c-Noetherian if and only if it is Noetherian, and it is c-principal if and only if it is principal.
    \item $\N$ is c-principal. In fact, fix a congruence $\cC$ and look at the collection of elements $S_\cC = \{(x,y)\in \N \times \N : y>0, x\sim x+y\}$. Then $\cC = \langle x_0 \sim x_0 + y_0\rangle$ where $x_0 = \min_{(x,y)\in S} x$ and $y_0 = \mathrm{gcd}(y: (x,y) \in S)$.
    \item As the proof above depends only on properties of addition, it follows that $\N_{\max} = \N \cup \{-\infty\}$ with $a\oplus b :=\max(a,b)$ and $a\otimes b :=a+b$ is also c-principal.
    \item $\N_{\max, \min} := \N \cup \{\pm \infty\}$ with $a\oplus b := \max(a,b)$ and $a\otimes b:= \min(a,b)$ is not c-Noetherian. In fact, for every $n\in \N$, $\{y\sim x : y, x \leq n\text{ or } y=x\}$ is a congruence. 
    \item Every finite semiring is trivially c-Noetherian, but they might not be c-principal.
    \item Let $\cM_b$ be the minmax semiring with $b$ elements, that is, the set $\{1,\ldots, b\}$ with $x+y = \min(x,y)$ and $xy=\max(x,y)$. Then $\cM_b$ is c-principal if and only if $b\leq 3$.
\end{enumerate}
\end{eg}

As expected, c-Noetherian semirings have equivalent formulations similar to Noetherian rings. This reconciles our concept with certain semirings discussed in Section 3.1 of \cite{JooMincheva2017}. The proof is almost the same as the ring theoretic analogue.

\begin{prop}
    \label{p:eqcnoeth}
    Let $R$ be a semiring. The following are equivalent.
    \begin{enumerate}
        \item $R$ is c-Noetherian.
        \item Congruences on $R$ satisfy the ascending chain condition (ACC). That is, any chain $C_1\subseteq C_2\subseteq \ldots \subseteq C_n \subseteq C_{n+1} \subseteq \ldots$ of congruences stabilizes (that is, $C_n =C_{n+1}$ for every $n\gg 1$).
        \item Any non-empty collection of congruences in $R$ has a maximal element. 
    \end{enumerate}
\end{prop}

\begin{proof}
    Assume (i). Let $\cC_1\subseteq \cC_2 \subseteq \ldots $ be an increasing chain. Let $\cC = \bigcup_i \cC_i$. It is immediate that $\cC$ is a congruence. As $R$ is c-Noetherian, $\cC$ is finitely generated, say $\cC = \langle a_i \sim b_i : i=1,\ldots, n\rangle$ for some $n\geq 1$. Then, for each $i$, there is some $n_i\geq 1$ such that $(a_i, b_i) \sim \cC_{n_i}$. Let $N=\max_i(n_i)$. Then $(a_i, b_i)\in \cC_N$ for each $i$, that is $\cC\subseteq \cC_N$, proving that $\cC_n = \cC$ for each $n\geq N$, proving (ii).

    Now assume (ii). Let $\cM$ be a non-empty collection of congruences on $R$. Suppose $\cM$ has no maximal elements. Then for any $\cC\in \cM$, we can find a congruence $\cC'$ strictly containing $\cC$. Thus, by Axiom of Dependent Choice, this means that we can find a strictly increasing chain of congruences $C_1\subsetneq C_2\subsetneq \ldots \subsetneq C_n \subsetneq C_{n+1} \subsetneq \ldots$ inside $\cM$, contradicting (ii). So $\cM$ must have a maximal element. 

    Finally, assume (iii). Let $\cC$ be a congruence and let $\mathcal{F}$ be the collection of all finitely generated congruences consisting of elements of $\cC$. By (iii), $\mathcal{F}$ has a maximal element, say $\cC'$. Now let $(a,b)\in \cC$. Then $\cC' \subseteq \langle \cC' \cup (a,b)\rangle \in \mathcal{F}$ implies, by the maximality of $\cC'$, that $\cC'=\langle\cC' \cup (a,b)\rangle$. Thus, $(a,b)\in \cC'$. This shows that $\cC \subseteq \cC'$, proving that $\cC = \cC'$ is finitely generated. 
\end{proof}

As in the world of rings, c-principal and c-Noetherian semirings are well-behaved under some standard semiring constructions. The following result summarises this fact.

\begin{prop}
    \label{p:quotprod}
    Let $A$ be a c-principal (resp. c-Noetherian) semiring.
    \begin{enumerate}
        \item For any c-principal (resp. c-Noetherian) semiring $B$, $A\times B$ is c-principal (resp. c-Noetherian).
        \item For any multiplicatively closed set $S$, $S^{-1}A$ is c-principal (resp. c-Noetherian).
        \item Any quotient $A \twoheadrightarrow B$, $B$ is c-principal (resp. c-Noetherian).
        \item Define the \textit{conical extension} $A_\omega$ as $A\cup \{\omega\}$ with the operations $a+\omega = a$ and $a\cdot \omega = \omega$ for each $a\in A$. Then every congruence of $A_\omega$ is generated by two (resp. finitely many) relations.
        \item $G(A) = \Z \otimes_\N A$ is a principal ring (resp. a Noetherian ring).
    \end{enumerate}
\end{prop}

\begin{proof}
    For (i), suppose $\cC$ is a congruence of $A\times B$. Let $\cC_A = \{(a,a')\in A^2: ((a,0), (a', 0))\in \cC\}$ and $\cC_B = \{(b,b')\in B^2: ((0,b), (0,b'))\in \cC\}$. It is immediate that $\cC_A, \cC_B$ are congruences on $A, B$ respectively, and that if $(a,a')\in \cC_A$ and $(b,b')\in \cC_B$ then $((a,b), (a',b'))\in \cC$. Now, suppose $((a,b),(a',b'))\in \cC$. Multiplying both sides with $(1,0)$ (resp. with $(0,1)$), it is immediate that $(a,a')\in \cC_A$ and $(b,b')\in \cC_B$. Thus, if $\cC_A =\langle (a_i, a'_i) : i\in I\rangle$ and $\cC_B =\langle (b_j, b'_j) : j\in J\rangle$, then $\cC = \langle (a_i, b_j) \sim (a'_i, b_j'):i\in I, j\in J\rangle$. The claim follows from here. 

    For (ii), we classify all congruences of $R =S^{-1}A$. Given a congruence $\cC$ of $A$, define $S^{-1}\cC$ as $\{ (a/s, a'/s')\in (S^{-1}A)^2 : \exists\> t\in S : (ts'a, tsa')\in \cC\}$. This is a congruence on $S^{-1}A$. Conversely, let $\cD$ be a congruence on $S^{-1}A$. Define $\cC$ as $\{(a,a')\in A^2 : (a/1, a'/1)\in \cD\}$. This is a congruence on $A$ and it can be seen that $S^{-1}\cC= \cD$. Thus, every congruence of $S^{-1}A$ comes from a congruence of $A$. Furthermore, if $\cC = \langle (a_i, a'_i) : i\in I\rangle$, then it is immediate that $\langle (a_i/1, a'_i/1) : i\in I\rangle$. This proves the claim.

    The proof of (iii) is immediate from the fact that congruences of $B$ correspond to the congruences of $A$ which contain the congruence kernel of the map. 

    For (iv). we claim that every congruence of $A_\omega$ looks like $\cC^{+} = \cC \cup \{(\omega, \omega)\} $ or $\cC^{++} =\cC \cup \{(\omega, \omega)\}\cup \{ (\omega , a), (a,\omega) : (a,0)\in \cC\}$ where $\cC$ is a congruence on $A$. Firstly, it is immediate that these are congruences. Conversely, let $\cC'$ be a congruence on $A_0$. It is clear that $\cC = \cC'\cap A^2$ is a congruence on $A$. If $\omega$ is only related to itself, then we are done. Otherwise, suppose $I = \{a\in A : (a, \omega)\in \cC\}$. Then $I$ is an ideal of $A$. In particular, $0\in I$. Furthermore, if $(a,0)\in \cC$ if and only if $a\in I$ by transitivity. Thus, $\cC' = \cC \cup \{(\omega, \omega)\}\cup \{ (\omega , a), (a,\omega) : (a,0)\in \cC\}$ as desired. Now suppose $\cC = \langle (a_i,b_i) : i\in I\rangle $ is a congruence on $A$. Then $\cC^{+} = \langle (a_i,b_i) : i\in I\rangle $ and $\cC^{++} = \langle \{(a_i,b_i) : i\in I\} \cup (0, \omega)\rangle$. This proves the claim.

    For (v), we claim that every ideal of $G(A)$ looks like $G(\cC) = \{1\otimes a + (-1)\otimes a' : (a,a')\in \cC\}$. Firstly, for a congruence $\cC$, $G(\cC)$ is indeed an ideal. Conversely, given an ideal $I$ of $A$, define $\cC(I)$ as $\{(x,y) \in A^2 : 1\otimes x +(-1)\otimes y \in I\}$. We claim that $\cC(I)$ is indeed a congruence on $A$. It is immediate it is an equivalence relation which is closed under addition. For multiplication, we observe that for any $(x,y), (x',y')\in \cC(I)$, we have
    $$1\otimes xx' + (-1)\otimes yy' = (1\otimes x)(1\otimes x' + (-1)\otimes y') + (1\otimes y')(1\otimes x + (-1)\otimes y)\in I$$
    proving that $(xx', yy')\in \cC(I)$. Finally, for any ideal $I$, one has that $G(\cC(I))= I$, proving that every ideal looks like $G(\cC)$. Finally, we want to prove that if $\cC = \langle (a_i, b_i) : i\in I\rangle$, then $G(\cC) = \langle 1\otimes a_i + (-1)\otimes b_i : i\in I\rangle $. That is, if $(a,b)\in \cC$, we wish to show that $1\otimes a + (-1)\otimes b \in \langle 1\otimes a_i + (-1)\otimes b_i : i\in I\rangle$. By Proposition \ref{p:definingconggen} and telescopic summation, it suffices to assume that $a\sim b$ is an elementary relation. So there are elements $u, v,c\in A$ and $i\in I$ such that $a = ua_i+vb_i+c$ and $b=ub_i+va_i+c$. This implies that
    $$1\otimes a + (-1)\otimes b = (1\otimes u + (-1)\otimes v)(1\otimes a_i + (-1)\otimes b_i)\in G(\cC).$$
    So $G(\cC) = \langle 1\otimes a_i + (-1)\otimes b_i : i\in I\rangle$. This proves our claim. 
\end{proof}

\begin{rem}
    In the process, we have proved the following well-known classification results. Let $A, B$ be two semiring. 
    \begin{enumerate}
        \item Any congruence of $A\times B$ appears uniquely as $\pi(\cC_A, \cC_B)$ where $\cC_A, \cC_B$ are congruences on $A, B$ respectively and $\pi : (A\times A) \times (B \times B) \to (A\times B)\times (A\times B)$ is the permutation map $((a,a'), (b,b'))\mapsto ((a,b),(a',b'))$.
        \item For a multiplicative subset of $A$, every congruence of $S^{-1}A$ looks like $S^{-1}\cC : =\{ (a/s, a'/s')\in (S^{-1}A)^2   : \exists\> t\in S : (ts'a, tsa')\in \cC\}$.
        \item For a congruence $\cD$ of $A$, the congruences of $A/\cD$ correspond one-to-one to the congruences of $A$ containing $\cD$.
        \item Any congruence of $A_\omega$ has the form $\cC^{+} = \cC \cup \{(\omega, \omega)\} $ or $\cC^{++} =\cC \cup \{(\omega, \omega)\}\cup \{ (\omega , a), (a,\omega) : (a,0)\in \cC\}$ where $\cC$ is a congruence on $A$.
        \item Any ideal of the ring $G(A) = \Z\otimes_\N A$ looks like $G(\cC) = \{1\otimes a +(-1)\otimes b : (a,b) \in \cC\}$ where $\cC$ is a congruence on $A$.
    \end{enumerate}
\end{rem}

\begin{rem}
     It is worth noting that even if $A$ is a c-principal semiring, it might happen that $A_\omega$ is not c-principal. For example, if $A=\N$ and $\cC = \langle 0 \sim \omega, 1\sim 2\rangle$, then $\cC$ is not principal. Observe that $A_\omega/\cC\simeq \bB$. Suppose $\cC$ is principal, say $\cC = \langle a\sim b \rangle$. We observe that as $\omega$ is only related to $0$ and itself, $0\sim \omega$ is an elementary relation with respect to $a\sim b$.
     
    So, there are elements $u,v,c\in B$ such that $ua+vb+c = \omega$ and $ub+va+c  =0$. But $B$ is negative free by definition. In particular, we must have $c=ua=vb=\omega$ and $0 = ub+va\neq \omega$. This can happen if and only if $a=\omega$ or $b=\omega$. Without loss of generality, assume that $a=\omega$. Then $b\sim \omega$ implies that $b=0$ or $b=\omega$. In the first case, the congruence has quotient $\N$ and in the latter, it is the diagonal congruence. This is a contradiction. Thus, $\cC$ is not principal.

    $A_\omega$ is c-principal if and only if $A$ is c-principal and each and that every congruence $\cC$ is a Rees congruence (see definition \ref{d:k-Noetherian}).
\end{rem}
    
The following proposition produces some new examples of c-Noetherian semirings. The proof only uses the monoid structure and has been known for a while. 

\begin{prop}
    \label{p:finc-Noetherian}
    Let $S$ be a $\N$-algebra that is finite as an $\N$-module. Then $S$ is c-Noetherian.
\end{prop}

\begin{proof}
    Suppose $S$ is not c-Noetherian, that is, there is an infinite chain $C_1\subsetneq C_2 \subsetneq C_3 \ldots$ of congruences on $S$. For a monoid $U$, let $\Z[U]$ be the free $\Z$-algebra generated by it. For example, if $U=\N$, then $\Z[U]\simeq \Z[T]$, the polynomial ring in one variable over $\Z$. Let $S_i = S/C_{i-1}$ for $i\geq 2$ and $S_1=S$ and $S_0$ be a free finitely generated $\N$-module surjecting on to $S_1$. Then the non-isomorphic surjections $S_0\twoheadrightarrow S_1 \twoheadrightarrow S_2\ldots$ yield non-isomorphic surjections $\Z[S_0]\twoheadrightarrow \Z[S_1] \twoheadrightarrow \Z[S_2]\twoheadrightarrow\ldots$ that do not terminate. But the first ring in this sequence is isomorphic to a finite polynomial ring over $\Z$, and is therefore Noetherian in the classical sense. This contradicts the existence of such a chain of surjections.
\end{proof}

At the time of writing, it was still not known whether any algebra over a c-Noetherian semiring $S$ which is finite as an $S$-module is c-Noetherian again. Finitely generated algebras are usually not c-Noetherian --- see Theorem \ref{thm1}.

It is natural to ask when a semiring being c-Noetherian implies that it is Noetherian, in the sense that every ideal is finitely generated. Unfortunately, there is ideals and congruences of an arbitrary semiring do not always recover each other. Instead, we look at a more restrictive class, the class of k-ideals, following \cite{abuhlail2019knoetheriankartiniansemirings}. We will see how they relate to the c-Noetherian condition later on.

\begin{defn}
    \label{d:k-Noetherian}
    Let $R$ be a semiring. We say that an ideal $I$ is a \textit{k-ideal} if $\mathrm{ker}(R\twoheadrightarrow R/\cC_I) = I$, where $\cC_I = \langle x\sim 0: x\in I\rangle$ is the Rees congruence of $I$. 

 In fact, $\cC_I = \{(x,y) : \text{there are }a,b\in I\text{ such that } x+a = y+b\}$. So an ideal $I$ is a k-ideal if and only if for each $x, y\in R$ with $x, x+y\in I$, we have $y\in I$.

    Given any subset $S\subseteq R$, the \textit{k-ideal generated by $S$} is defined to be the smallest k-ideal containing $S$. We say $R$ is \textit{k-Noetherian} if every k-ideal is finitely generated (as a k-ideal).
 \end{defn}

 \begin{rem}
     \label{r:nomkideal} Even though the term k-ideal and k-Noetherian have already appeared in the literature, there have been multiple interpretations of what k means. Of many possible interpretations, we feel the most appropriate expansion is either \textit{kernel} (as such ideals appear as kernels) or \textit{k\"undbar} (which means cancellable in German).
 \end{rem}

The following proposition has a proof almost identical to Proposition \ref{p:eqcnoeth} and hence the proof is omitted.

 \begin{prop}
     \label{p:k-Noeth}
     Let $R$ be a semiring.
     \begin{enumerate}
         \item The k-ideal generated by an ideal $I$ is the kernel of the map $R\twoheadrightarrow R/\cC_I$.
         \item The following are equivalent.
         \begin{enumerate}
             \item $R$ is k-Noetherian.
             \item The k-ideals satisfy the ascending chain condition.
             \item Any non-empty collection of k-ideals has a maximal element.
         \end{enumerate}
     \end{enumerate}
 \end{prop}

 The following proposition says that k-ideals behave somewhat better than arbitrary ideals in semirings of a special kind.

 \begin{prop}
    \label{p:k-ideal}
    Let $R$ be a totally ordered ring and $S=R_{\geq0}$. Then there is a bijective correspondence
    $$\{\text{Ideals of }R\} \longleftrightarrow \{\text{k-ideals of }S\}$$
    given by 
    $$I\subseteq R \mapsto I_{\geq 0}\subseteq S$$
    and
    $$J\subseteq S \mapsto JR = \{x-y : x,y\in J\} \subseteq R.$$   
    Finally, if $I$ is generated by $\cD$, then $I_{\geq 0}$ is generated, as a k-ideal, by $\{|x| :  x\in \cD\}$. Thus, if $R$ is Noetherian, then $S$ is k-Noetherian.
\end{prop}

\begin{proof}
    These maps are inverses of each other and map to the correct domains. Now, if $I = R\langle S\rangle$ and $x\in \cD$. Then $|x|\in I_{\geq 0}$. Now, let $J$ be the k-ideal generated by $|\cD| = \{|x| : x\in \cD\}$. Then $I:=JR= R\langle|\cD|\rangle = R\langle \cD\rangle$, proving that $I_{\geq 0} = J$ and the claim about generation. k-Noetherianness follows immediately.
\end{proof}

\section{A Non-c-Noetherian Semiring}\label{s4}

In this section, we explicitly construct a somewhat universal example of a Non-c-Noetherian semiring.

\begin{thm}\label{thmNc-Noetherian}
    Let $R= \bB[X]$. Then the congruence $\cC=\langle X^n+1\sim X^m+1 : n, m\geq 1\rangle$ is not finitely generated. 
\end{thm}

\begin{proof}
    The proof is split into two observations. Firstly, we observe that for each $p\geq 0$, and for any $a\in R$, $X^p\sim a$ in $\cC$ implies that $a=X^p$. Indeed, by Proposition \ref{p:definingconggen}, it suffices to assume that $a\sim X^p$ is an elementary relation with respect to the given generating set. So, there are $n,m\geq 1$, $u,v, c\in R$ such that $X^p = u(X^n+1)+v(X^m+1) + c$ and $a = v(X^n+1)+u(X^m+1)+c$. However, observe that if $x,y\in R$ are such that $x+y=X^p$, then $x,y$ are either $X^p$ or $0$. As $X^t+1 \nmid X^p$ for $t\geq 1$, we must have $u(X^n+1), v(X^m+1)=0$, which implies that $u=v=0$. This shows that $X^p = c = a$, as desired.
    
    Now, suppose, for the sake of contradiction, that $\cC$ was finitely generated. Say $\cC=\langle a_1\sim b_1,\ldots,a_k\sim b_k\rangle$ with each $a_i\neq b_i$. In particular, $a_i\neq 0, X^p$ for any $p\geq 0$. Let $N$ be greater than the degrees of the finitely many polynomials $a_i,b_i$.

    In a spirit similar to the first paragraph, we claim that if $X^N+1\sim a$, then $a=X^N+1$. Once again, we may assume that $a$ is related to $X^{N}+1$ by an elementary relation with respect to the given generating set $\{(a_i,b_i) : i=1,\ldots, k\}$. So there is an integer $1\leq i\leq k$ and $u,v,c\in R$ such that 
    $X^N+1 = ua_i+vb_i + c$ and $a = va_i+ub_i+c$.

    Now any summand of $X^N+1$ is one of $0,1, X^N, X^N+1$. In particular, $ua_i$ is one of these. But, as $X^N+1$ is irreducible, and $a_i$ does not divide $1, X^N$ or $X^{N}+1$ due to degree constraints and the first paragraph, we must have $ua_i=0$. Thus $u=0$. Similarly, we can show that $v=0$. Thus, $X^N+1 = c =a$. 
    
    This shows that if $\cC$ is finitely generated, then $X^N+1 \sim a$ implies $X^N+1 = a$. But $X^N+1 \sim X^{N+1}+1$ in $\cC$. This is a contradiction to our hypothesis. Thus, $\cC$ is not finitely generated. 
\end{proof}

\begin{cor}
     The semiring $\bB[X]/\cC$ with $\cC$ as above is finitely generated as an algebra over $\N$ but not finitely presented.
\end{cor}

\begin{proof}
    Immediate from Theorem \ref{thmNc-Noetherian} and Proposition \ref{p3}.
\end{proof}

\begin{rem}
    It must be noted that $\bB[X]$ is not even Noetherian in the sense that there are ideals of $\bB[X]$ that cannot be finitely generated. For instance, consider the ideal generated by $\{X^{n}+1 : n\geq 1\}$. An argument similar to the one given above proves that this ideal is not finitely generated. However, it is k-Noetherian. Indeed, let $I$ be a nonzero k-ideal. Let $n\geq 0$ be the least such that $I$ has an element of the form $f(X) = X^n +(\text{higher degree terms})$. Then $X^n+f(X) = f(X)\in I \implies X^n\in I$. Also, as for any other term in $I$, every monomial has degree at least $n$, we get that $I=(X^n)$, showing that every k-ideal is actually principal (as an ideal, not just a k-ideal).
\end{rem}

\begin{rem}
    In \cite{JooMincheva2018} Proposition 2.12, it is shown that $\bB[X,Y]$ and $\bB(X,Y)$ are not c-Noetherian. This appears to be a predecessor to the notion of c-Noetherian semirings. However, it turns out $\bB(X)$ is, in fact, c-Noetherian.
\end{rem}

We now prove Theorem \ref{thm1}.

\thilb*

\begin{proof}
    The implications $(iv)\implies (iii) \implies (ii) \implies (i)$ are immediate. Now suppose $(i)$ holds. For the sake of contradiction, suppose that $R$ is not a ring. Then $R$ admits $\bB$ as a homomorphic image by Corollary \ref{thBJ}. Then we have a surjective map $\phi : R[X] \to \bB[X] \to \bB[X]/\cC$, where $\cC$ is the congruence in Theorem \ref{thmNc-Noetherian}. But then the congruence kernel $\cC_\phi$ is not finitely generated as $\cC$ is not finitely generated. This contradicts our assumption that $R[X]$ is c-Noetherian. So, $R$ must be a ring for which $R[X]$ is Noetherian. Thus, $R = R[X]/(X)$ is a Noetherian ring.
\end{proof}

\section{Positive Models and Finiteness}\label{s7}

One of the major motivations to study semirings come from the study of positive models. While many different notions exist in the literature, here is what we will mean by a positive model construction. For alternate views with the intent of geometric applications, the reader is directed to \cite{borger_jun_semiring}, \cite{FockGoncharov2006}, and \cite{ToenVaquie2009}.

\begin{defn}\label{d:positivity function}
    Let $R$ be a ring. A \textit{positivity} or \textit{trichotomy function} is a function $p: R\to \{0, \pm 1\}$ such that
    \begin{enumerate}
        \item $p(0)=0, p(1)=1$.
        \item $p(xy)=p(x)p(y)$.
        \item $p(x+y)=-1\implies p(x)=-1$ or $p(y)=-1$.
    \end{enumerate}

    In such a case, the set $R_p = \{x\in R : p(x)\geq 0\}$ is a subsemiring of $R$.
\end{defn}

The most common examples of trichotomy functions are the signum function $\text{sgn}:\R \to \{0,\pm 1\}$ and the zero divisor function $p(x) = \mathbf{1}\{x \text{ is a nonzerodivisor}\}.$ 

\begin{defn}\label{d:positive model}
    Let $R$ be a ring endowed with a positivity function $p$ and let $S=R_p$. Let $\cC$ be a subcategory of $\mathbf{Mod}_R$. A \textit{positive model construction} on $\cC$ is a functor $F : \cC \to \mathbf{Mod}_S$ with a natural isomorphism $ R \otimes_{S} F \xrightarrow{\sim} Id$.
\end{defn}

The following is the most important construction for a positive model for the sake of our considerations.

\begin{eg}
    Let $R=\Z$, $p=\text{sgn}$. Let $\cC$ be the category of all $\Z$-subalgebras of $\R$, with the morphisms being sign-preserving ring homomorphisms. Define $F(R)=R_{\geq 0}$. Then one may easily verify that $F$ is a positive model construction. As no choice of a basis is involved, this positive model is called the \textit{canonical positive model} of a subring of $\R$. Similar constructions can be given for any totally ordered ring. 
\end{eg}

\begin{eg}\label{eg:tensorposp-1}
    Suppose $R$ is a ring endowed with a positivity function such that $p(-1)=-1$. Then we claim that for any $R_p$-module $M$, we have $R\otimes_{R_p} M \simeq \Z \otimes_\N M$.

    By transitivity of base change, it suffices to prove that $R \simeq R':=\Z\otimes_\N R_p$. Given any $x\in R$, either $p(x)\geq 0$ or $p(-x)=-p(x)\geq 0$. So the map $\pi : R' \to R$ given by $n\otimes r \mapsto nr$ is indeed a surjection. Now define the set-theoretic map $\phi: R \to R'$ given by $\phi(x) = 1\otimes x$ if $p(x)\geq 0$ and $(-1)\otimes (-x)$ if $p(x)<0$. It is clear that $\pi \circ \phi = Id_R$. Given any $r\in \Z\otimes_\N R_p$, by linearity, $r = 1\otimes x+(-1)\otimes y$ for some $x,y\in R_p$. Now let $u =\pi(r) = x-y$. If $p(u)\geq 0$, then 
    $$r =  1\otimes (u+y) + (-1)\otimes y = 1\otimes u + (1+(-1))\otimes y = 1\otimes u.$$ 
    If $p(u) = -1$, then a similar argument shows $r = (-1)\otimes (-u)$. In either case, $r$ is in the image of $\phi$. Thus $\phi \circ \pi \circ \phi (x) = \phi(x)$ implies that $\phi \circ \pi = Id_{R'}$.
    
    In fact, it follows that the isomorphism $R\otimes_{R_p} M \simeq \Z \otimes_\N M$ is functorial in $M$.
\end{eg}

\begin{eg}
    Another example of a positivity structure is as follows. Let $R$ be a ring equipped with a homomorphism $f: R \to \R$. Define the \textit{induced positivity structure} on $f$ as $p_f(x) = \mathrm{sgn}(f(x))$. It is immediate that $p$ is a positivity structure. For example, let $a\in \Z$ and $ev_a : \Z[X] \to \Z$ be the evaluation map. Then we can look at $p(x) := \mathrm{sgn} \circ ev_a$. A similar construction can be done for any ring $S$ with an existing positivity structure in place of $\R$.
\end{eg}

We have the following rigidity result on trichotomy functions.

\begin{prop}
    \label{p:trichotomyimpliessign}
    Suppose $p$ is a positivity structure on a ring $R$. Then $p(-1)=1$ if and only if $R_p=R$. In particular, if $R$ has positive characteristic, then any positivity function on $R$ is non-negative. 
\end{prop}

\begin{proof}
    If a ring $R$ has positive characteristic $n>1$, then $p(-1)=p(n-1)\geq 0$ by property (iii) of a positivity function. Thus, the second conclusion is immediate from the first one.

    It is also clear that if $p(-1)=\pm 1$, and that $p(-1)=-1$ implies $R_p\neq R$. So assume $p(-1)=1$ and $R_p\neq R$. Choose an $r\in R$ such that $p(r) = -1$. As $p(r) = p((r+1)+(-1)) = -1$, we must have $p(r+1)=-1$. So $p(r^2+r) = p(r+1)p(r) = 1 = p(r^2)$. Thus, $r^2+r, r^2, -1 \in R_p$ and thus $r=(r^2+r)-r^2 \in R_p$. This is a contradiction. So, if $p(-1)=1$, then we must have $R_p=R$.
\end{proof}

\begin{prop}\label{p:posmod}
    Let $R$ be a ring equipped with a positivity function $p$ and $S = R_p$. Given any $S$-module $M$, $R\otimes_S M$ is finitely generated/finitely presented/free/projective/flat as an $R$-module if so is $M$ as an $S$-module. Furthermore, if $M$ is an $S$-algebra, then $M\otimes_S R$ is a Noetherian/principal ring if $M$ is c-Noetherian/c-principal over $S$.
\end{prop}

\begin{proof}
    The first part is immediate as finitely generated/finitely presented/free/ projective/flat modules are preserved by base change. Now, if $M$ is an $S$-algebra, we may assume that $p(-1)=-1$ by Proposition \ref{p:trichotomyimpliessign}. By the computations in Example \ref{eg:tensorposp-1}, it suffices to show that $\Z\otimes_\N M$ has the desired properties. But this is part (v) of Proposition \ref{p:quotprod}.
\end{proof}

In particular, if a positive model has a good property, then we can expect the ambient object to have the analogous nice property. The converse, however, is far from true. For example, $\Z/\N$ is definitely not flat over $\N$, but $\Z\otimes_\N \Z \simeq \Z$ is flat over $\Z$.

Our main result for the rest of the article is the following, which shows that the converse is true in a broad class of rings.

\ttarget*

\begin{rem}
    We claim that Theorem \ref{t:target} is, in a sense, the best that one could hope from canonical positive models. Let $R\neq \Z$ be a subring of $\R$ and $S=R_{\geq 0}$. We claim that $S$ is never finitely generated or free as a module over $\N$. The finite generation is immediate --- any finitely generated $\N$-submodule of $\R_{\geq 0}$ is necessarily discrete, while $S$ is always dense in $\R_{\geq 0}$. For freeness, suppose $S$ is free over $\N$ with basis $\mathcal{B}$. Then $\mathcal{B}$ has at least two elements $u, v$. As $u\neq v$, we can assume, without loss of generality, that $u>v$. Then $u-v\in S\setminus\{0\}$ implies that there is a $k\geq 1$, vectors $v_1,\ldots, v_k\in \mathcal{B}$ and natural numbers $\lambda_1,\ldots, \lambda_k>0$ such that $u-v = \sum_{k} \lambda_kv_k$. So $u = \sum_{k} \lambda_k v_k + v$, contradicting that $\mathcal{B}$ is a basis.
\end{rem}

For the rest of this section, we work towards proving the claims related to finite generation of congruences. 

\begin{thm}\label{t:fin}
    Suppose $R\subseteq \R$ is a subring and $S=R_{\geq 0}$. Let $\sim$ be a congruence on $S$. Given any $x\in S$, define $I_x = \{y\in S : y+x\sim x\}$. Then
    \begin{enumerate}
        \item $I_x$ is a k-ideal.
        \item If $x\leq y$, then $I_x\subseteq I_y$.
        \item If $I_t\neq 0$, then $I_x=I_t$ for each $x\geq t$.
        \item If $R\neq \Z$, then $I_x=I_y$ for each $x,y>0$.
    \end{enumerate}
\end{thm}

\begin{proof}

For (i), we fix am $x\in S$. Clearly, $0\in I_x$. Also, if $a,b\in I_x$, then
    $$x\sim x+b \sim (x+a)+b$$
    shows that $a+b\in I_x$. Similarly $a,b\in S$ and $a, a+b\in I_x$ implies that
    $$x+b \sim (x+a)+b \sim x+(a+b) \sim x,$$
   Finally, let $a\in I_x, a\neq 0$ and $r\in S$. Let $N\in \N$ be large enough so that $Na>rx$. Then 
    \begin{align*}
        x+ra&\sim (x+Na)+ra \\
        &= (x+Na-rx)+r(x+a)\\
        &\sim (x+Na-rx)+rx \\
        &\sim x+Na \sim x
    \end{align*}
    proving that $I_x$ is indeed an ideal. 
    
    For (ii), suppose $x\leq y$. Now $a\in I_x$ implies that 
    $$y = (y-x)+x \sim (y-x)+(x+a) \sim y+a$$
    proving that $a\in I_y$. So $I_x\subseteq I_y$.

    For (iii), suppose that $I_t\neq 0$. Let $y\in I_t\setminus 0$. Let $x\geq t$ and $a\in I_x$. Choose $N\geq 1$ large enough so that $Ny>x$. Then 
    \begin{align*}
   t+a &\sim (t+Ny)+a \\
   &\sim (t+ (a+x)) + (-x+Ny) \\
   &\sim (t+x) + (-x+Ny) \\
   &=t+Ny \sim t.
    \end{align*}
   This shows that $I_x=I_t$, proving our claim.

   Finally, for (iv), if $R\neq \Z$, then $R$ is dense in $\R$. Let $x>y>0$. Then there is a $\beta\in S, \beta >0$ such that $\beta x<y$. Then $\beta I_x \subseteq I_{\beta x} \subseteq I_y \subseteq I_x$. So $I_x=0$ if and only if $I_y=0$. Furthermore, (iii) shows that $I_x=I_y$.
\end{proof}

This theorem has a few interesting consequences that we will now proceed to prove. Recall that the conical extension $A_\omega$ from Proposition \ref{p:quotprod}. Observe that we always have a quotient map $\pi : A_\omega \to A$ given by $\omega \mapsto 0$ and $a\mapsto a$ for each $a\in A$. Furthermore, if $A$ does not have zero divisors or negatives, then we may define a section $s: A \to A_{\omega}$ given by $a\mapsto a$ if $a\neq 0$ and $0 \mapsto \omega$. 

\begin{rem}[Rees-Type Congruences]
    Let $R$ be a totally ordered ring and $S = R_{\geq 0}$. Given any ideal $I$ of $R$, define the \textit{Rees congruence} $\cC_I$ on $S$ as $\{(x,y) \in S^2:  x-y\in I\}$ and the \textit{altered Rees congruence} $\cC_I^{-} := \{(x,y)\in S^2 :x,y>0, x-y\in I\}\cup \{(0, 0)\}$. It is immediate that these are congruences. Furthermore, if $I=0$, then $\cC_I=\cC^{-}_I = \Delta$, the diagonal congruence. In all other cases, that is, when $I\neq 0$, $\cC_I\supsetneq \cC^{-}_I$, hence the name. 

    Observe that $S/\cC_I \simeq R/I$ and $S/\cC^{-}_I \simeq (R/I)_\omega$.
\end{rem}

While these are congruences in the positive model of any totally ordered ring, they gain additional importance for subrings of $\R$..

\begin{cor}
    \label{c1:t:fin}
    Suppose that $R, S$ is as in Theorem \ref{t:fin} and $R\neq \Z$. Then the only congruences on $S$ are of the form $\cC_I$ or $\cC^{-}_I$ for an ideal $I$.
    Furthermore, if $I_{\geq 0}$ is generated as a k-ideal (see Proposition \ref{p:k-ideal}) by $\{x_i:i\in K\}$, then $\cC_I$ is generated by $\langle 0 \sim x_i : i\in K\rangle$ and $\cC^{-}_I$. is generated by $\langle 1 \sim 1+x_i : i\in K\rangle$
\end{cor}

\begin{proof}
    From Theorem \ref{t:fin} and the remark above, it is clear that a congruence is classified by the information of whether $I_0=0$ or $I_0\neq 0$, and what $I_1$ is. Let $I=I_1R$ If $I_0 = 0$, send then the congruence is simply given by $\cC^{-}_I$. If $I_0\neq 0$, then the congruence is $\cC_I$ by Theorem \ref{t:fin}.

    Now for the generation, the claim is immediate for $\cC_I$. For $\cC^{-}_I$, we observe that the congruence $\cD = \langle 1 \sim 1+x_i : i\in K\rangle$ is contained in $\cC^{-}_I$ and they have identical $I_1$ and $I_0$. Thus, by Theorem \ref{t:fin}, they must agree. 
\end{proof}

\begin{cor}
    \label{c2:t:fin}
    Suppose $R$ is a subring of $\R$ and $S=R_{\geq 0}$. Then $S$ is c-Noetherian if and only if $S$ is k-Noetherian, which holds if and only if $R$ is Noetherian. Furthermore, if $R$ is a PID, then $S$ is c-principal. Finally, if $R$ is a real order, then any nontrivial quotient of $S$ is finite. 
\end{cor}

\begin{proof}
    The claims are immediate for $R=\Z$, so we assume $R\neq \Z$. Observe that, by Corollary \ref{c1:t:fin}, $S$ is c-Noetherian if and only $S$ is k-Noetherian, which, by Proposition \ref{p:k-ideal} and Proposition \ref{p:posmod}, is true if and only if $R$ is Noetherian. As this is indeed the case, we are done. Now if $R$ is a PID, then the last line of Corollary \ref{c1:t:fin} implies that $S$ is c-principal. Finally, if $R$ is a real order, then the result on quotients of $S$ in Corollary \ref{c1:t:fin} implies that every nontrivial quotient of $R$ is finite.
\end{proof}

\begin{rem}
    What is also immediate from basic algebraic number theory is that if $R$ is a ring of integers (maximal order) to start with, then every congruence on $S$ is generated by at most two elements and that every non-trivial quotient is principal. 
    
    Also, for $R=\Z[\pi]$, we find that $S=R_{\geq 0}$ is c-Noetherian, unlike $S'=\N[\pi]\simeq N[T]$. This gives us some more evidence that perhaps the canonical positive model is better than a basis-dependent model.
\end{rem}

This settles the finiteness part of Theorem \ref{t:target}. However, we point out that this is not the only approach or proof of the theorem. While discussing a complete alternative proof leads to an unnecessary detour, we state a few lemmas, sketch the proof and make a remark to justify its inclusion. Fix $R\neq \Z$.

\begin{lem}
    \label{lem7.1}
    Let $A$ be a semiring, $\sim$ a congruence on $A$, and suppose that $x\sim x+y$. Then for each $n\geq 1$, we have $x^n \sim x^n +y^n$.
\end{lem}

\begin{proof}
    We proceed by induction. The base case is assumed. Suppose that the claim is true for $n\geq 1$. Then 
    $$x^{n+1} \sim x(x^n+y^n) \sim x^{n+1} + y^n(x+y)\sim x(x^n+y^n)+y^{n+1}\sim x^{n+1}+y^{n+1}.$$
    The claim is thus established for each $n\geq 1$ by induction.
\end{proof}

\begin{lem}
    \label{lem7.2}
       Let $R\subseteq \R$ be any algebraic $\Z$-algebra and $\sim$ be a non-trivial congruence in $S=R_{\geq 0}$. Then there are integers $n>m>1$ such that $n\sim m$.
\end{lem}

\begin{proof}
    Let $a \sim a+u$ in $S$ with $a,u>0$. Then, by Lemma \ref{lem7.1}, for each $n\geq 1$, we have  $a^n \sim a^n + u^n$. So for any $q(T)\in T\N[T]$, $q(a)\sim q(a)+q(u)$. Now $u$ satisfies an equation of the form $uf(u) = ug(u) + l$ where $f,g\in \N[X], l\in \N\setminus \{0\}$. So,
    \begin{align*}
        af(a)+ag(a)+l   &\sim af(a)+ag(a) +ug(u)+l \\
                        &=af(a)+ag(a) +uf(u) \\
                        &\sim af(a)+ag(a).
    \end{align*}

Now, let $m$ be any integer exceeding $af(a)+ag(a)$ and $n=m+l$. As $\sim$ is a congruence, we have 
\begin{align*}
m &= (m-af(a)-ag(a)) + af(a)+ag(a) \\
&\sim (m-af(a)+ag(a)) + af(a)+ag(a) + l \\
&= m+l =n.    
\end{align*}
\end{proof}

We now sketch the alternate proof of the theorem for real orders.

\begin{enumerate}
    \item As $R \neq \Z$ is a real order, we can find for $k=\dim_\Z R\geq 2$ elements $\alpha_1,\ldots, \alpha_k\in R, \alpha_i\geq 0$ such that $R = \Z\alpha_1 \oplus \ldots \oplus \Z\alpha_k$.
    \item Now, for every element $x\in S, x>0$, choose $y\in S, y\neq 0$ such that $my<x$. Write $y=\sum_i \lambda_i \alpha_i$. Assume, without loss of generality, that $\lambda_1\geq 1$.
    \item Now using $m\sim n$ repeatedly, we can show that $my \sim m\sum_i \mu_i \alpha_i$ where each $\mu_i\geq 1$.
    \item Finally, by Euclid's division algorithm and some minor adjustments, one ends up showing that $x\sim \sum_i \theta_i \alpha_i$, where each $\theta_i$ is a positive integer smaller than $n$, proving the finiteness. The proof of c-principal is similar.
\end{enumerate}

\begin{rem}
    It is tempting to try to generalize the above results further. While we are using only certain specific properties of real numbers, we rely on it being totally ordered and Archimedean. However, it follows that a totally ordered Archimedean ring $R$ is isomorphic to a subring of the real numbers. Indeed, consider the map $\phi : R \to \R$ as $x\mapsto \sup \{ r/s :  r\cdot 1_R \leq s\cdot 1, r\in \Z, s\in \N\}$. This is easily seen to be an injective ring homomorphism. So, in a sense, subrings of $\R$ are the best cases in which the above techniques work. The same is also valid for cancellative totally ordered Archimedean semirings - these are canonical positive models for some subring of real numbers. On the other hand, Lemma \ref{lem7.1} and a weaker version of Lemma \ref{lem7.2} can be proved without these strong assumptions, giving us an alternative approach and justifying their inclusion.
\end{rem}

\section{Flatness of a Positive Model of Real Orders}\label{s8}

We now prove the following result to complete the proof of Theorem \ref{t:target}. 

\begin{thm}
    For any subgroup $X\subseteq \R$, $X_{\geq 0}$ is strongly flat over $\N$.
\end{thm}

The only proofs that we are aware of at this point involve some convex methods and analysis. The proof essentially consists of two parts, namely

\begin{enumerate}[label=(P\arabic*)]
    \item Translation of the problem into a problem about semilattices in $\R^n$, and
    \item Homothetic projection onto a hyperplane and its properties.
\end{enumerate}

Let us do (P1) now. We start with the following lemma.

\begin{lem}
    \label{lem8.1}
    To prove the theorem, it suffices to show that for any $X$ finitely generated as a module over $\Z$, $X_{\geq 0}$ is strongly flat over $\N$.
\end{lem}

\begin{proof}
    $X$ is a directed union of its finitely generated submodules $\{X_\alpha\}_\alpha$. It is immediate that $X_{\geq 0}$ is the directed union of the collection $\{(X_{\alpha})_{\geq 0}\}$. If each $(X_\alpha)_{\geq 0}$ is strongly flat, then so is $X_{\geq 0}$. This proves our claim. 
\end{proof}

Henceforth, we assume that $X$ is finitely generated, hence free. Let $\{\alpha_1,\ldots, \alpha_k\}$ be a $\Z$-basis of $X$ such that $\alpha_i> 0$ for each $i$. Then $X_{\geq 0}$ can be identified with the semilattice 
    $S=\{(n_1,\ldots, n_k)\in \Z^k : \sum_i n_i\alpha_i\geq 0\}$. 
    
Define $\gamma := (\alpha_1,\ldots,\alpha_n)\in \R^n$ and the hyperplane $L := \{\mathbf{v} \in \R^n : \mathbf{v}\cdot \gamma = 0\}$. Then $S = \{\mathbf{v} \in \Z^n : v\cdot \gamma >0\} \cup \{0\}$. We show that every finitely generated submodule $M$ of $S$ is contained in a finitely generated free submodule of $S$.

Let $\mathbf{V}=\{\mathbf{v} \in \R^n : v\cdot \gamma >0\} \cup \{0\}$ and let $S'$ be the set of elements of $S$ with coprime coordinates. 

\begin{lem}\label{lem8.2}
   Suppose $M\subseteq S$ is finitely generated as a module over $\N$. To show that $M$ is contained in a finitely generated free submodule of $S$, it is sufficient to find vectors $v_1,\ldots, v_n\in \Z^n$ such that
    \begin{enumerate}
        \item $(v_1|v_2|\ldots|v_n)\in GL(n,\Z)$,
        \item $v_i\cdot \gamma >0$, and
        \item $\R_{\geq 0}\langle v_1,\ldots,v_n\rangle \supseteq M$.
    \end{enumerate}
\end{lem}

\begin{proof}
    The first two conditions shows that $v_i \in S$ and that $\N\langle v_1,\ldots,v_n\rangle$ is free by Lemma \ref{lem8.3}. Finally, if $m\in M$, then 
    \begin{align*}
        &m\in \R_{\geq 0}\langle v_1,\ldots,v_n\rangle \\
        &\implies  m = \sum_{i=1}^{n} \lambda_i v_i \text{ for some }\lambda_i\in \R_{\geq 0} \\
        &\implies (\lambda_1,\ldots,\lambda_n) =\sum_{i=1}^{n} \lambda_ie_i = X^{-1}\sum_{i=1}^{n} \lambda_iv_i = X^{-1}m \in \Z^n.
    \end{align*}
So $\lambda_i\in \Z\cap \R_{\geq 0} = \N$. Hence, $M\subseteq \N\langle v_1,\ldots,v_n\rangle$, proving our claim. 
\end{proof}

\begin{lem}
    \label{lem8.3}
    Let $S$ be a cancellative semiring, $R=S\otimes_\N\Z$, $M$ an $S$-module, and $M' = N\otimes_S R$. If $M$ is generated by $n$ elements and $M'\simeq R^n$, then $M$ is generated freely by these $n$ generators.
\end{lem}

\begin{proof}
    Consider the diagram
    \[\begin{tikzcd}
        S^n \arrow[d, two heads, "f"]\arrow[r, hook] & R^n \arrow[d, two heads, "f_R"] \arrow[dr, "g"] & \\
        M \arrow[r] & M' \arrow[r, "\sim"] & R^n
    \end{tikzcd}\]
    As $g:R^n \to R^n$ is surjective and $R^n$ is a free $R$-module, $g$ is an isomorphism. Hence, $f_R$ is an isomorphism, which immediately implies that $f$ is injective, proving our claim.
\end{proof}

For simplicity, we make the following definition.

\begin{defn}
    A collection of vectors $V=\{v_1,\ldots, v_n\}$ in $\Z^n$ is called \textit{$\gamma$-oriented} if 
    \begin{enumerate}
        \item $(v_1|\ldots|v_n) \in GL_n(\Z)$, and
        \item $v_i\cdot \gamma >0$.
    \end{enumerate}
    
    For $R=\R_{\geq 0}, \N$, its $R$-span is denoted by $\mathrm{Span}_R(V)$. Given a $\gamma$-oriented collection $V_0 = \{v_1,\ldots,v_n\}$, an \textit{elementary refinement} is either of the following:
    \begin{enumerate}
        \item For distinct indices $1\leq i \neq j \leq n$ such that $v_i\cdot \gamma > v_j\cdot \gamma$, a collection $V'=\{v_1',\ldots, v_n'\}$ of the form $v'_k =v_k $ if $k\neq i$ and $v'_i = v_i - v_j$.
        \item A permutation $V^\sigma_0 = \{v_{\sigma(1)},\ldots, v_{\sigma(n)}\}$.
        \end{enumerate}
    
    We say that another $\gamma$-oriented collection $W=\{w_1,\ldots,w_n\}$ \textit{refines} $V$ if there is a finite sequence of $\gamma$-oriented collections $\{V_n\}_{n=0}^{N}$ such that each $V_i$ is an elementary refinement of $V_{i-1}$, $V_0=V$ and $V_N=W$. We will denote this as $V\preceq W$.
\end{defn}

\addtocounter{defn}{-1}
\addtocounter{lem}{3}

It is immediate that $\preceq$ is a partial preorder and that $V\preceq W \preceq V$ if and only if $V$ is a permutation of the indices of $W$.

\begin{lem}\label{lem3}
    To prove the theorem, it suffices to show that given any $\gamma$-oriented collection $W$ and a vector $m\in S$, there is another $\gamma$-oriented collection $V$ refining $W$ such that $m\in \mathrm{Span}_{\R_{\geq 0}}(V)$.
\end{lem}

\begin{proof}
    This follows from Lemma \ref{lem8.2} and an easy induction.
\end{proof}

\begin{rem}
    It is tempting to think that $V \preceq W$ if and only if $\mathrm{Span}_\N(V) \subseteq \mathrm{Span}_\N(W)$. Although such a claim is true for $\Z$-spans, or even $\N$-spans when $n=\dim_
    \Z(X) = 2$, this claim is not true for $n\geq 3$. In fact, consider $V = \{ (5,0,3), (1,2,1), (1,1,0)\}$ and let $W = \{(1,0,0), (0,1,0), (0,0,1)\}$. Then we have $\mathrm{Span}_{\N}(V) \subseteq \mathrm{Span}_\N(V_s)$ but $V\npreceq W$. This example easily lifts to give examples for each $n\geq 3$.
\end{rem}

The rest of this section is dedicated to proving the hypothesis of the preceding lemma. That is, given any $\gamma$-oriented collection $W$ and a vector $m\in S$, we find another $\gamma$-oriented collection $V$ such that $\R_{\geq 0}\langle v_1,\ldots,v_n\rangle \supseteq \{w_1,\ldots, w_n,m\}$. 

To do this, we move on to part (P2). To study $\R_{\geq 0}$-spans, it suffices to study their homothetic projections onto a fixed hyperplane. So let 
\begin{align*}
    H = &\{x\in \R^n : x\cdot \gamma =1\}, \\
    p: &\mathbf{V} \longrightarrow H \\
    & x \longmapsto \frac{x}{x\cdot\gamma}
\end{align*}
So, for $v, w\in S$, $p(v)=p(w)$ if and only if $\mu v=\lambda w$ for some $\mu, \lambda\geq 1$. In particular, if $v,w\in S'$, then $p(v)=p(w)$ if and only if $v=w$. Hence, we may define an inverse map $g : p(S') \to S'$. Also, for $v\in V$, define $|v|_\gamma = v\cdot \gamma$. If $v\in \Z^n$, then $|v|_\gamma = 0 \iff v=0$ as $\alpha_1,\ldots, \alpha_n$ are linearly independent over $\Z$. Thus, $|.|_\gamma$ can be used to define a total order on $\Z^n$.

Given a $\gamma$-oriented collection $W=\{w_1,\ldots,w_n\}$, define $\cM(W)$ as the set of all refinements of $W$. By the previous lemma, it suffices to show that the union

\begin{equation}\label{eq:union}\tag{6.1}
    \bigcup_{\{v_1,\ldots,v_n\}\in \cM(W) } \R_{\geq 0}\langle v_1,\ldots,v_n\rangle = \mathbf{V}.
\end{equation}

By taking projections, this is equivalent to showing
    $$\cO(W) := \bigcup_{\{v_1,\ldots,v_n\}} p(\R_{\geq 0}\langle v_1,\ldots,v_n\rangle) = \bigcup_{\{v_1,\ldots,v_n\}}\cv (p(v_1),\ldots,p(v_n)) = H, $$
where $\cv(.)$ is the convex hull and the union is taken over $\cM(W).$

To prove that Equation (\ref{eq:union}) holds, we prove multiple lemmas.

\begin{lem}
    \label{lem4}
    Given any $\gamma$-oriented collection $W=\{w_1,\ldots,w_n\}$ and a collection of indices $J\subseteq \{1,\ldots, n\}$ with $|J|\geq 2$ and any $\delta>0$, there is a $\gamma$-oriented collection $\{v_1,\ldots, v_n\}\in \cM(W)$ such that
    \begin{enumerate}
        \item $v_i=w_i$ if $i\notin J$,
        \item $|v_i|_\gamma < \delta $ for every $i\in J$, and
        \item $w_i \in \N\langle v_j | j\in J\rangle$ for each $i\in J$.
    \end{enumerate}
\end{lem}

\begin{proof}
    Fix two distinct indices $j_1,j_2 \in J$ and let $x_i = w_{j_i}$. Without loss of generality, let $|x_1|_\gamma < |x_2|_\gamma $ Let $G=\{M\in GL(2,\Z) | M^{-1}\in M(2,\N)\}$. Then $G$ is closed under multiplication. We will now produce inductively a sequence of vectors $x^{(n)} = \{x_1^{(n)}, x_2^{(n)}\}$ such that
    \begin{enumerate}
        \item $x^{(0)}_i = x_i$ for $i=1,2$,
        \item $|x^{(n)}_2|_\gamma>|x_1^{(n)}|_\gamma>0$ for $n\geq 0$,
        \item $x^{(n+1)}$ is a refinement of $x^{(n)}$, and
        \item If $m_n = x^{(n)}_2 \cdot \gamma$, then $m_{n+1} < m_{n}/2$.
    \end{enumerate}
The base case is immediate. So we proceed by induction. Assume that we are given $x^{(n)}$. Let $m\geq 1$ be such that $m|x^{(n)}_1|_\gamma<|x_2^{(n)}|_\gamma<(m+1)|x_1^{(n)}|_{\gamma}$. The injectivity of $x \mapsto x\cdot \gamma$ on $\Z^n$ ensures a strict inequality. Let $x^{(n+1)}_2 = x^{(n)}_2-mx_1^{(n)}$. Clearly $0<|x^{(n+1)}_2|_\gamma < |x_1^{(n)}|_\gamma$. So, let $m'\geq 1$ be such that
$m'|x^{(n+1)}_2|_\gamma<|x_1^{(n)}|_\gamma<(m'+1)|x_2^{(n+1)}|_{\gamma}$ and let $x_1^{(n+1)} = x_1^{(n)} - mx_2^{(n+1)}$. It follows that $0<|x_1^{(n+1)}|_\gamma< |x_2^{(n+1)}|_\gamma$ and that $x^{(n+1)}$ is a refinement of $x^{(n)}$. Finally, $2m_{n+1} < (1+mm')|x_2^{(n+1)}|_\gamma + m |x_1^{(n+1)}|_\gamma = |x_2^{(n)}|_\gamma = m_n$, as desired.

By choosing $n$ large enough, we have found a matrix $M = \begin{pmatrix}
    a& b\\c & d
\end{pmatrix}$ in $G$ and vectors $v_1, v_2$ such that $0<|ax_1+bx_2|_\gamma,|cx_1+dx_2|_\gamma<\delta$. Let $v_{j_1}=ax_1+bx_2$, $v_{j_2}=cx_1+dx_2$. For $j\in J$, $j\neq j_i$, choose $k_j$ such that $k_j |v_{j_1}|_\gamma< |w_j|_\gamma<(k_j+1) |v_{j_1}|_\gamma$ and let $v_j = w_j - k_jv_1$. Also, for $j\notin J$, let $v_j=w_j$. It is immediate that this collection has the desired properties.
\end{proof}

\begin{lem}
    \label{lem5}
    Given any two indices $1\leq i < j \leq n$, the entire line through the points $p(w_i)$ and $p(w_j)$ is in $\cO(W)$.
\end{lem}

\begin{proof}
    By Lemma \ref{lem4}, we can construct a sequence $W^{(n)} = \{w_k^{(n)}\}$ of $\gamma$-oriented collections such that $w^{(n+1)}\in \cM(w^{(n)})$, $w_k^{(n)} = w_k$ if $k\neq i, j$ and $0<w_k^{(n)}\cdot \gamma <\frac{1}{n}$. Let $x_n=w^{(n)}_i, y_n = w^{(n)}_j$. From the construction in Lemma \ref{lem4}, it is clear that there are integers $a_n, b_n, c_n, d_n\geq 0$ such that $a_nd_n-b_nc_n=1, x_n = a_nx_{n+1}+b_ny_{n+1}$ and $y_n = c_nx_{n+1}+d_ny_{n+1}$. Thus, the line segments $L_n$ that join $p(x_n)$ and $p(y_n)$ form an ascending chain. Thus, it suffices to show that the usual absolute values $|p(x_n)|, |p(y_n)|$ diverge to infinity. But 
    $$|p(x_n)| = \frac{|x_n|}{x_n\cdot \gamma} > n|x_n| \geq n$$ as $x_n$ has integer coordinates. So we are done. 
\end{proof}

Note that, in fact, this completes the proof when $n=2$. We point out a notable similarity of the arguments presented in the previous two lemmas to the theory of continued fractions and convergents.

\begin{lem}
    \label{lem6}
    Let $n\geq 3$. Given any $V = \{v_1,\ldots,v_n\} \in \cM(W)$ and any pair of indices $1\leq i,j\leq n$, such that $p(w_i)>p(w_j)$ there is a $U=\{u_1,\ldots,u_n\} \in \cM(V)$ such that $p(w_i-w_j) \in \cv(p(u_1),\ldots,p(u_n)).$
\end{lem}

\begin{proof}
    Without loss of generality, assume that $i=1, j=2$. As $V$ is a $\Z$-basis of $\Z^n$, there are integers $\mu_k$ such that $e =w_1-w_2 = \sum_k \mu_k v_k$. Let $\gamma_k=v_k\cdot \gamma >0$. So, 
    $$\sum_k \mu_k \gamma_k = e\cdot \gamma>0.$$
    Suppose $Q$ is the number of $1\leq k\leq n$ such that $\mu_k<0$.

\underline{\textbf{CASE I: $Q=0$}} \\
Just choose $U=V$.

\underline{\textbf{CASE II: $Q=1$}} \\
By possibly permuting the vectors, assume that $\mu_1 = -s<0$. Let $\delta = \frac{1}{sn}\min\{e\cdot \gamma, \gamma_1 \}>0$. By Lemma \ref{lem4}, we get a $U'=\{u_1',\ldots,u_n'\}\in \cM(V)$ such that $0<|u_j'|_\gamma<\delta$ for each $j\geq  2$ and $u_1'=v_1$. Furthermore, $v_k \in \mathrm{Span}_\N\langle u_j'|j\geq 2\rangle$ for $k\geq 2$. So, $w_1-w_2 =\sum_j \mu'_j
u_j'$ where $\mu_1'=\mu_1$ and $\mu_j'\geq 0$ if $j\geq 2$.

Consider the set of vectors $(\lambda_2,\ldots, \lambda_n)$ of non-negative integers such that 
\begin{enumerate}
    \item $\lambda_js \leq \mu'_j$ if $j\geq 2$
    \item If $\eta = \sum_{j\geq 2} \lambda_j u_j'$, then $\gamma_1>\eta\cdot \gamma$, and
\end{enumerate}

The set of such vectors $(\lambda_2, \ldots, \lambda_n)$ is finite and non-empty (as the zero vector satisfies these conditions). Choose a termwise maximal element of this set. Thus, for each $j\geq 2$, we must have either $\lambda_j = \left\lfloor\frac{\mu_j'}{s}\right\rfloor =:\lambda_j^{max}$ or $(\eta+u_j')\cdot \gamma > \gamma_1$. Otherwise, the vector $\lambda' = \lambda+e_{j}$, where $e_j$ is the $j$-th standard basis vector, satisfies all the three conditions, contradicting maximality.

Let $J_1 = \{2\leq j \leq n : \lambda_j=\lambda_j^{max}\}$ and $J_2 = \{2\leq j \leq n : j\notin J_1\}$.

Observe that 
\begin{align*}
\sum_{j\geq 2} \lambda_j^{max}u_j'\cdot \gamma
    &> \sum_{j\geq 2} \frac{\mu_j'u_j'}{s}\cdot \gamma - u_j'\cdot \gamma \\
    &= \frac{e\cdot \gamma}{s} + \gamma_1 - \sum_j u_j'\cdot \gamma \\
    &> \frac{e\cdot \gamma}{s} + \gamma_1 - n\delta > \gamma_1.
\end{align*}

So $J_2\neq \emptyset$. Also, for any $j\in J_2$, we have $\mu'_j \geq (\lambda_j+1)s$.

Define
$$u_j = \begin{cases}
    u_1' -\eta & j=1 \\
    u_j' & j\in J_1\\
    u_j'+\eta - u_1' & j\in J_2 \\
\end{cases}$$

It is immediate that $U =\{u_1,\ldots, u_n\} \in \cM(U')\subseteq \cM(V)$. Furthermore, 
$$w_1 - w_2 = \left(-s+\sum_{j\in J_2} (\mu'_j  - s\lambda_j)\right)  u_1 +\sum_{j\geq 2} (\mu'_j  - s\lambda_j)u_j$$
is in $\mathrm{Span}_\N(U)$. This proves our claim.

\underline{\textbf{CASE III: $Q\geq 2$}} \\
After possibly permuting the indices, suppose $\mu_1,\ldots, \mu_t<0$ and $\mu_{t+1},\ldots, \mu_n\geq 0$. Let $M=\sum_{k>t}\mu_k>0$ (as $e\cdot \gamma>0)$ and $\delta = {M}^{-1}(e\cdot \gamma)>0$. By Lemma \ref{lem4}, choose $U'=\{u_1',\ldots,u_n'\}\in \cM(V)$ such that $0<|u_j'|_\gamma<\delta$ for each $j\leq t$ and $u_j'=v_j$ if $j>t$. So, $e = \sum_k \mu_k'u'_k$, where $\mu_k'<0$ if $k\leq t$ and $\mu_k'=\mu_k$ for $k>t$.

We claim the existence of a $Z =\{z_1,\ldots,z_n\} \in \cM(U')$ such that if 
$e = \sum \theta_k z_k$, then $\theta_k=\mu_k'$ if $k\neq t$ and $\theta_t > 0$. Rest will follow from induction and a reduction to Case II.

For $k>t$, let $\lambda_k$ be the largest integers such that $(u_k' - \lambda_k u_t')\cdot \gamma >0$. Then
\begin{align*}
    \mu'_t+\sum_{k>t} \mu'_k\lambda_k &\geq \mu'_t+\sum_{k>t} \mu'_k\left(\frac{u_k'\cdot \gamma}{u_t'\cdot \gamma}-1\right) \\
    &= \frac{\sum_{k<t} (-\mu'_k)u'_k\cdot \gamma + e\cdot \gamma }{u'_t\cdot \gamma} - M \\
    &\geq \frac{e\cdot \gamma}{\delta } - M > 0 
\end{align*}

Now let 
$$z_i = \begin{cases}
    u'_i & i\leq t \\
    u'_i - \lambda_iu_t' & i>t.
\end{cases} $$

Then $Z\in \cM(U')$ and we have
\begin{align*}
    e &= \sum_k \mu_k'u'_k \\
    &= \sum_{k<t} \mu_k'u_k'+ \sum_{k>t} \mu_k'(u_k'- \lambda_k u_t') + \left(\mu'_t+\sum_{k> t}\mu_k\lambda_k\right) u_t'. \\
    & = \sum_{k\neq t} \mu_k' z_k + \theta_t z_t
\end{align*}

where $\theta_t = \sum_{k>t}\mu_k\lambda_k + \mu_t' > 0$. 
\end{proof}

\begin{lem}
\label{lem7}
    For any $\gamma$-oriented $W$, $\cM(W)$ is a directed set under refinement.
\end{lem}

\begin{proof}
    We show that given any refinement $V$ of $W$ and an elementary refinement $W'$ of $W$, there is a common refinement of $V, W'$. If $W'$ is a permutation, then $V$ itself is the common refinement. Otherwise, we may use Lemma \ref{lem6} directly.
\end{proof}
    
\begin{lem}
\label{lem8}
    For any $\gamma$-oriented $W=\{w_1,\ldots,w_n\}$, $\cO(W)$ is $H$.
\end{lem}

\begin{proof}
    From Lemma \ref{lem7}, $\cM(W)$ is a directed collection and hence $\{p(V) :  V \in \cM(W)\}$ is a directed collection of convex sets. So their union, $\cO(W)$ is a convex set as well. Additionally, from Lemmas \ref{lem5} and \ref{lem6}, we get that $\cO(W)$ contains $n-1$ independent lines (in both directions) starting at any fixed vector $v$. The only such convex set is, indeed, the whole of $H$, proving our claim.
\end{proof}

This completes the proof of the theorem. 

\begin{rem}
    A very similar proof will show that whenever $S, S'$ are subrings of $\R$ with $S\subseteq S'$, then $S'_{\geq 0}$ is strongly flat over $S_{\geq 0}$.
\end{rem}

{\footnotesize \textbf{Acknowledgments:} The author thanks Professor James Borger, his PhD supervisor, for his guidance and support. Most of the problems and the work presented here were consequences of several fruitful discussions on semirings and their theory with him. The author also thanks the Australian National University and the committee of Deakin PhD Scholarship for supporting his PhD, both academically and financially, as the work was entirely done during the course.}

\bibliographystyle{plain}
\bibliography{fgsr}

@article{borger_jun_semiring,
  author       = {Borger, James and Jun, Jaiung},
  title        = {Facets of Module Theory over Semirings},
  journal      = {Advances in Mathematics},
  volume       = {479},
  pages        = {110431},
  year         = {2025},
  eprint       = {2405.18645},
  archivePrefix= {arXiv},
  primaryClass = {math.RA}
}

@misc{abuhlail2019knoetheriankartiniansemirings,
      title={On k-{N}oetherian and k-{A}rtinian Semirings}, 
      author={Jawad Abuhlail and Rangga Ganzar Noegraha},
      year={2019},
      eprint={1907.06149},
      archivePrefix={arXiv},
      primaryClass={math.RA},
      url={https://arxiv.org/abs/1907.06149}, 
}

@article {allencohen,
    AUTHOR = {Allen, Paul J.},
     TITLE = {Cohen's theorem for a class of {N}oetherian semirings},
   JOURNAL = {Publ. Math. Debrecen},
  FJOURNAL = {Publicationes Mathematicae Debrecen},
    VOLUME = {17},
      YEAR = {1970},
     PAGES = {169--171},
      ISSN = {0033-3883,2064-2849},
   MRCLASS = {16A78},
  MRNUMBER = {302709},
MRREVIEWER = {R.\ Gorton},
       DOI = {10.5486/pmd.1970.17.1-4.19},
       URL = {https://doi.org/10.5486/pmd.1970.17.1-4.19},
}

@article {katsov2004,
    AUTHOR = {Katsov, Yefim},
     TITLE = {On flat semimodules over semirings},
   JOURNAL = {Algebra Universalis},
  FJOURNAL = {Algebra Universalis},
    VOLUME = {51},
      YEAR = {2004},
    NUMBER = {2-3},
     PAGES = {287--299},
      ISSN = {0002-5240,1420-8911},
   MRCLASS = {16Y60 (06A12 18A40 18G05)},
  MRNUMBER = {2099802},
MRREVIEWER = {Jan\ Paseka},
       DOI = {10.1007/s00012-004-1865-1},
       URL = {https://doi.org/10.1007/s00012-004-1865-1},
}

@article {allenideal2006,
    AUTHOR = {Allen, Paul J. and Neggers, Joseph and Kim, Hee Sik},
     TITLE = {Ideal theory in commutative {$A$}-semirings},
   JOURNAL = {Kyungpook Math. J.},
  FJOURNAL = {Kyungpook Mathematical Journal},
    VOLUME = {46},
      YEAR = {2006},
    NUMBER = {2},
     PAGES = {261--271},
      ISSN = {1225-6951,0454-8124},
   MRCLASS = {16Y60},
  MRNUMBER = {2235566},
MRREVIEWER = {Sam\ L.\ Blyumin},
}

@incollection {Borger2016,
    AUTHOR = {Borger, James},
     TITLE = {Witt vectors, semirings, and total positivity},
 BOOKTITLE = {Absolute arithmetic and {$\Bbb F_{1}$}-geometry},
     PAGES = {273--329},
 PUBLISHER = {Eur. Math. Soc., Z\"urich},
      YEAR = {2016},
      ISBN = {978-3-03719-157-6},
   MRCLASS = {13F35 (11T30 14A20)},
  MRNUMBER = {3701903},
MRREVIEWER = {Ramdorai\ Sujatha},
}

@article {allenhilb1980,
    AUTHOR = {Dale, Louis and Allen, Paul J.},
     TITLE = {An extension of the {H}ilbert basis theorem},
   JOURNAL = {Publ. Math. Debrecen},
  FJOURNAL = {Publicationes Mathematicae Debrecen},
    VOLUME = {27},
      YEAR = {1980},
    NUMBER = {1-2},
     PAGES = {31--34},
      ISSN = {0033-3883,2064-2849},
   MRCLASS = {16A78},
  MRNUMBER = {592180},
MRREVIEWER = {G\"unter\ R.\ Krause},
       DOI = {10.5486/pmd.1980.27.1-2.06},
       URL = {https://doi.org/10.5486/pmd.1980.27.1-2.06},
}

@article{Borger_2015,
   title={Boolean {W}itt vectors and an integral {E}drei–{T}homa theorem},
   volume={22},
   ISSN={1420-9020},
   url={http://dx.doi.org/10.1007/s00029-015-0198-6},
   DOI={10.1007/s00029-015-0198-6},
   number={2},
   journal={Selecta Mathematica},
   publisher={Springer Science and Business Media LLC},
   author={Borger, James and Grinberg, Darij},
   year={2015},
   month=nov, pages={595–629} }

@article{JooMincheva2017,
  author    = {Dániel Joó and Kalina Mincheva},
  title     = {Prime congruences of idempotent semirings and a Nullstellensatz for tropical polynomials},
  journal   = {Journal of the London Mathematical Society},
  volume    = {96},
  number    = {3},
  pages     = {--},
  year      = {2017},
  eprint    = {1408.3817},
  archivePrefix = {arXiv},
  primaryClass   = {math.AC},
  note      = {Introduces prime congruences via twisted products for idempotent semirings},
}

@article{JooMincheva2018,
  author    = {Dániel Joó and Kalina Mincheva},
  title     = {On the dimension of polynomial semirings},
  journal   = {Journal of Algebra},
  volume    = {507},
  pages     = {103--119},
  year      = {2018},
  doi       = {10.1016/j.jalgebra.2018.04.007},
  note      = {Focuses on prime congruences and Krull dimension for additively idempotent semirings},
}

@phdthesis{Zelich2026,
  author = {Ivan Diklich-Zelich},
  title = {The Miracle of Flatness in Algebraic Geometry},
  school = {Columbia University},
  year = {2026},
  eprint = {2607.17406},
  archivePrefix = {arXiv},
  primaryClass = {math.AG}
}

@book{FreeseMcKenzieMcNultyTaylor2022,
  author    = {Ralph S. Freese and Ralph N. McKenzie and
               George F. McNulty and Walter F. Taylor},
  title     = {Algebras, Lattices, Varieties, Volume II},
  series    = {Mathematical Surveys and Monographs},
  volume    = {268},
  publisher = {American Mathematical Society},
  address   = {Providence, RI},
  year      = {2022},
  isbn      = {978-1-4704-6797-5}
}

@misc{chen2026arithmeticsemiringsiideals,
      title={The Arithmetic of Semirings Part I: Ideals}, 
      author={Jay Chen and Trevor Hyde and Dorien Laurens and Jasper Piermarini and Harriet Simons},
      year={2026},
      eprint={2607.10778},
      archivePrefix={arXiv},
      primaryClass={math.AC},
      url={https://arxiv.org/abs/2607.10778}, 
}

@article{FockGoncharov2006,
  author  = {Vladimir Fock and Alexander Goncharov},
  title   = {Moduli spaces of local systems and higher {Teichm\"uller} theory},
  journal = {Publications Math{\'e}matiques de l'IH{\'E}S},
  volume  = {103},
  year    = {2006},
  pages   = {1--211},
  doi     = {10.1007/s10240-006-0039-4}
}

@book{MaclaganSturmfels2015,
  author    = {Diane Maclagan and Bernd Sturmfels},
  title     = {Introduction to Tropical Geometry},
  series    = {Graduate Studies in Mathematics},
  volume    = {161},
  publisher = {American Mathematical Society},
  address   = {Providence, RI},
  year      = {2015}
}

@article{ToenVaquie2009,
  author  = {Bertrand To{\"e}n and Michel Vaqui{\'e}},
  title   = {Au-dessous de {Spec Z}},
  journal = {Journal of K-Theory},
  volume  = {3},
  number  = {3},
  year    = {2009},
  pages   = {437--500},
  doi     = {10.1017/is008004027jkt048}
}

@misc{GualdiKuhrsMayoGarciaXarles2026,
  author        = {Roberto Gualdi and Arne Kuhrs and Mayo Mayo Garc{\'i}a
                   and Xavier Xarles},
  title         = {Scheme theory for commutative semirings},
  year          = {2026},
  eprint        = {2601.14136},
  archivePrefix = {arXiv},
  primaryClass  = {math.AG}
}

@incollection{ConnesConsani2011,
  author    = {Alain Connes and Caterina Consani},
  title     = {Characteristic 1, entropy and the absolute point},
  booktitle = {Noncommutative Geometry, Arithmetic, and Related Topics},
  publisher = {Johns Hopkins University Press},
  address   = {Baltimore, MD},
  year      = {2011},
  pages     = {75--139}
}

@article{Lorscheid2017,
  author  = {Oliver Lorscheid},
  title   = {Blue schemes, semiring schemes, and relative schemes after
             {To{\"e}n} and {Vaqui{\'e}}},
  journal = {Journal of Algebra},
  volume  = {482},
  year    = {2017},
  pages   = {264--302},
  doi     = {10.1016/j.jalgebra.2017.03.023}
}

\end{document}